\documentclass[12pt]{article}

\title{On Continuity Properties of the Law \\
of Integrals of L\'evy Processes} 
\author{Jean Bertoin\thanks{Laboratoire de Probabilit\'es, 
Universit\'e Paris VI, 175 rue du Chevaleret, 75013 Paris,
France, email:
jbe@ccr.jussieu.fr} \and Alexander Lindner\thanks{ Centre for
Mathematical
Sciences, Technical University of Munich, Boltzmannstra{\ss}e 3,
D-85747
Garching, Germany, email: lindner@ma.tum.de} 
 \and Ross Maller\thanks{Centre for Mathematical Analysis,
 and School of Finance \& Applied Statistics, 
 Australian National University, Canberra, ACT,
 email: Ross.Maller@anu.edu.au} }

\usepackage{amssymb}
\usepackage{amsfonts}
\usepackage{amsmath}
\usepackage{amsthm}
\usepackage{graphics}
\usepackage{latexsym}

\newtheorem{theorem}{Theorem}[section]
\newtheorem{lemma}[theorem]{Lemma}

\newtheorem{example}[theorem]{Example}
\newtheorem{proposition}[theorem]{Proposition}

\newcommand{\halmos}{\quad\hfill\mbox{$\Box$}}

\def\topr{\buildrel P \over \to }

\def\=dr{\buildrel D \over = }

\def\nexto{\kern -0.54em}

\def\DoubleR{{\rm {I\ \nexto R}}}

\def\DoubleT{{\cal T \kern -0.7em T}}
\def\DoubleZ{{\cal Z \kern -0.7em Z}}
\def\DoubleC{{\rm\hbox{C \kern-0.8em\raise0.2ex\hbox{\vrule
height5.4pt width0.7pt}\kern0.2em}}}
\def\DoubleQ{{\rm\hbox{Q \kern-0.85em\raise0.25ex\hbox{\vrule
height5.4pt width0.7pt} \kern0.2em}}}

\def\LL{L\'{e}vy }

\newcommand{\supp}{\mbox{\rm supp}\,}

\begin{document}

\maketitle

\begin{abstract}
Let $(\xi,\eta)$ be a bivariate L\'evy process such that
the integral $\int_0^\infty e^{-\xi_{t-}} \, d\eta_t$ converges
almost
surely. We characterise, in terms of their \LL measures, those
L\'evy
processes  for which (the distribution of) this integral has
atoms. We then
turn attention to almost surely convergent integrals of the form
$I:=\int_0^\infty g(\xi_t) \, dt$, where $g$ is a deterministic
function. We
give sufficient conditions ensuring that $I$ 
has no atoms, and under further conditions derive that $I$
has a Lebesgue density. The results are also extended
to certain integrals of the form $\int_0^\infty g(\xi_t) \,
dY_t$, where $Y$
is an almost surely strictly increasing stochastic process,
independent of
$\xi$. 
\end{abstract}

\section{Introduction} \label{s1}
The aim of this paper is to study continuity properties of
stationary 
distributions of generalised Ornstein-Uhlenbeck processes and of
distributions of random variables of the form $\int_0^\infty g(\xi_t) \, dt$ for a
L\'evy process
$\xi$ and a general function $g: \mathbb{R}\to \mathbb{R}$.

For a bivariate L\'evy process $(\xi,\zeta) = (\xi_t,\zeta_t)_{t
\geq 0}$, the
{\it generalised Ornstein-Uhlenbeck (O-U) process} $(V_t)_{t
\geq 0}$ is defined as
$$V_t = e^{-\xi_t} \left( \int_0^t e^{\xi_{s-}} \, d\zeta_s + V_0
\right), \quad 
t \geq 0,$$
where $V_0$ is a finite random variable, independent of
$(\xi,\zeta)$.
This process appears as a natural continuous time generalisation
of random recurrence
equations, as shown by de~Haan and Karandikar \cite{HK}, and has
applications in many areas, such 
as risk theory (e.g. Paulsen~\cite{paulsen1}), perpetuities 
(e.g. Dufresne~\cite{dufresne}), financial time series (e.g.
Kl\"uppelberg
et al.~\cite{klm1}) or option pricing 
(e.g. Yor~\cite{yor1}), to name just a few.  See also Carmona et
al.~\cite{cpy1, cpy2}
for further properties of this process. Lindner and
Maller~\cite{lm} have shown
that the existence of a stationary solution to the generalised
O-U process is closely
related to the almost sure convergence of the  stochastic
integral
$\int_0^t e^{-\xi_{s-}} \, d\eta_s$ as $t\to \infty$, where
$(\xi,\eta)$ is a
bivariate L\'evy process, and $\eta$ can be explicitly
constructed in terms of 
$(\xi,\zeta)$. The stationary distribution is then given by
$\int_0^\infty e^{-\xi_{s-}}
\, d\eta_s$. Necessary and sufficient conditions for the
convergence of 
$\int_0^\infty e^{-\xi_{s-}} \, d\eta_s$ were obtained by
Erickson and Maller~\cite{em1}.
Distributional properties of the limit variable and hence of the
stationary distribution
of generalised O-U processes are of particular interest.
Gjessing and 
Paulsen~\cite{gp}  determined the distribution in many cases
when 
$\xi$ and $\eta$ are independent and the L\'evy measure of
$(\xi,\eta)$ is finite.
Carmona et al.~\cite{cpy1} considered the case when $\eta_t=t$
and the jump part of
$\xi$ is of finite
variation. Under some additional assumptions, they showed that
$\int_0^\infty 
e^{-\xi_{s-}} ds$ is absolutely continuous, and its density
satisfies 
a certain integro-differential equation.
In Section~\ref{s2} we shall be concerned with 
continuity properties of the 
limit variable $\int_0^\infty e^{-\xi_{s-}} \, d\eta_s$
without any restrictions on $(\xi,\eta)$, assuming only 
convergence of the integral. We shall give  
a complete characterisation of when this integral has atoms,
in terms of the characteristic
triplet of $(\xi,\eta)$.
This  characterisation relies on a similar result of 
Grincevi\v{c}ius~\cite{grin1} for ``perpetuities"
which are a kind of discrete time analogue of L\'evy integrals.

Then, in Section~\ref{s3}, 
we turn our attention to  continuity properties of the distribution of
the integral  $\int_0^\infty g(\xi_t) \, dt$, 
where
$\xi = (\xi_t)_{t \geq 0}$ is a one-dimensional L\'evy process
with non-zero L\'evy measure 
and $g$ is a general deterministic Borel function.
Such integrals appear in a variety of situations, for example concerningshattering phenomena in fragmentation processes, see, e.g.,Haas~\cite{haas}.
Fourier analysis and Malliavin calculus are classical tools for
establishing the absolute continuity of distributions of 
functionals of stochastic processes. In a different direction, 
the book of Davydov et al.~\cite{dls} treats
three different methods for proving absolute continuity of such
functionals: the ``stratification method'', the ``superstructure
method'' and the ``method of differential operators''. Chapter 4
in \cite{dls}
pays particular attention to Poisson functionals, which includes
integrals of L\'evy processes. While it may be hard to 
check the conditions and apply 
these methods in general (in particular to find admissible
semigroups for the stratification method), it has been carried
out in some cases.
For example, Davydov~\cite{davydov} gives sufficient conditions
for absolute 
continuity of integrals of the form $\int_0^1 g(X_t)\, dt$ for
{\it strictly
stationary} processes $(X_t)_{t \geq 0}$ and quite general $g$.
Concerning integrals of L\'evy processes,
Lifshits~\cite{lifshits}, p.~757, has shown that $\int_0^1
g(\xi_t) \, dt$
is absolutely continuous if $\xi$ is a L\'evy process with
infinite
and absolutely continuous L\'evy measure, 
and $g$ is locally Lipschitz-continuous and such that on a set
of 
full Lebesgue measure  in $[0,1]$ the derivative
$g'$ of $g$ exists and is continuous and non-vanishing;
see also 
Problem~15.1 in \cite{dls}. For our study of atoms of the
distributions
of integrals such as $\int_0^\infty g(\xi_t) \, dt$,
we will impose less restrictive assumptions on $g$ in
Section~\ref{s3}. 
Note also
that \cite{dls} and the references given there are usually
concerned with the
absolute continuity of functionals such as 
$\int_0^1 g(\xi_t) \, dt$ on the {\it compact} interval $[0,1]$,
while we
are concerned with integrals over $(0,\infty)$.
That absolute continuity 
of the distribution of integrals over compact sets and over
$(0,\infty)$ 
can be rather different topics 
is straightforward by considering the special case of compound Poisson processes.
See also part (iii) of Theorem~\ref{thm1} below for situations
where the integral over every finite time horizon may be
absolutely continuous, while the limit variable can degenerate to a constant.

Section~\ref{s3} is organised as follows: we start with some
motivating examples, in some of which
$\int_0^\infty g(\xi_t) \, dt$ has atoms while in others it
does not.  Then,
in Section~\ref{s3-2} we present some general criteria which
ensure the 
continuity of the distribution of $\int_0^\infty g(\xi_t) \,
dt$. The proofs there are
based on the sample path behaviour and on 
excursion theory for L\'evy processes. Then, in
Section~\ref{s3-3} we use
a simple form of the stratification method to obtain absolute 
continuity of $\int_0^\infty g(\xi_t) \, dt$ for certain cases
of $g$ and $\xi$ (which assume however no differentiability
properties
of $g$); the results are also extended to more general integrals
of
the form $\int_0^\infty g(\xi_t) \, dY_t$, where $Y=(Y_t)_{t
\geq 0}$
is a strictly increasing stochastic process, 
independent of the L\'evy process
$\xi$. 

Observe that
our focus  will be on continuity properties of the distribution
of
the integral  $\int_0^\infty g(\xi_t) \, dt$ (or similar
integrals), under
the assumption that it is finite a.s. A highly
relevant question
is to ask under which conditions the integral does converge. It
is important
of course that any conditions we impose to ensure continuity of
the
integral, or its absence, be compatible with  convergence. We
only
occasionally address this issue, when it is possible to give
some simple
sufficient (or, sometimes, necessary) conditions for
convergence. Our
approach is essentially to assume convergence and study the
properties of
the resulting integral. For a much fuller discussion of
conditions for
convergence {\it per se} we refer to
Erickson and Maller~\cite{em1, em2}, who give an overview of
known results
as well as new results on the finiteness of \LL integrals.

We end this section by setting some notation.
Recall that a {\it L\'evy process} $X=(X_t)_{t \geq 0}$ in
$\mathbb{R}^d$ ($d \in 
\mathbb{N}$) is a stochastically continuous process having
independent
and stationary increments, which has almost surely c\`adl\`ag
paths and 
satisfies $X_0 = 0$. For each L\'evy process, there exists a
unique 
constant $\gamma = \gamma_X = (\gamma_1, \ldots, \gamma_d) \in
\mathbb{R}^d$, 
a symmetric 
positive semidefinite matrix $\Sigma = \Sigma_X$, and a L\'evy
measure $\Pi =\Pi_X$
on $\mathbb{R}^d \setminus \{ 0 \}$ satisfying
$\int_{\mathbb{R}^d} 
\min \{ 1, |x|^2 \} \, \Pi_X(dx) < \infty$, such that for all $t
> 0$ 
and $\theta \in \mathbb{R}^d$ we have
$$(1/t) \log E \exp (i \langle \theta, X_t \rangle )= 
 i \langle \gamma, \theta \rangle - \frac{1}{2} \langle \theta,
\Sigma \theta \rangle
+ \int_{\mathbb{R}^d} (e^{i \langle z, \theta \rangle} -1 - i
\langle z, \theta
\rangle \mathbf{1}_{|z| \leq 1} )\, \Pi_X(dz ).$$
Here, $\langle \cdot, \cdot \rangle$ and $| \cdot |$ denote the
inner product
and Euclidian norm in $\mathbb{R}^d$, and $\mathbf{1}_{A}$ is
the indicator function of
a set $A$. Together, 
$(\gamma,\Sigma,\Pi)$ form the {\it characteristic triplet} of
$X$.
 The Brownian motion part of $X$ is described by the covariance
matrix $\Sigma_X$.
If $d=1$, then we will also write $\sigma^2_X$ for $\Sigma_X$,
and if $d=2$ and
$X=(\xi,\eta)$,  the upper and lower diagonal elements of
$\Sigma_X$ 
are given by $\sigma_\xi^2$ and $\sigma_\eta^2$. 
We refer to Bertoin~\cite{bertoin} and Sato~\cite{sato} for
further definitions
and basic properties of L\'evy processes. 
Integrals of the form $\int_a^{b} e^{-\xi_{t-}} \, d \eta_t$ for
a bivariate
L\'evy process $(\xi,\eta)$ are interpreted as the usual
stochastic integral
with respect to its completed natural filtration as in
Protter~\cite{protter},
where $\int_a^b$ denotes integrals over the set $[a,b]$, and
$\int_{a+}^b$ denotes
integrals over the set $(a,b]$. If $\eta$ (or a more general
stochastic process $Y
= (Y_t)_{t \geq 0}$ as an integrator) is of bounded variation on
compacts, then the stochastic
integral is equal to the pathwise computed Lebesgue-Stieltjes
integral, and will
also be interpreted in this sense. Integrals such as
$\int_0^\infty$ are to be
interpreted as limits of integrals of the form
$\int_0^t$ as $t\to\infty$, where the convergence will typically
be  almost sure.
The {\it jump} of a c\`adl\`ag process $(Z_t)_{t \geq 0}$ at
time $t$
will 
be denoted by $\Delta Z_t := Z_t - Z_{t-} = Z_t - \lim_{u
\uparrow t} Z_u$,
with the convention $Z_{0-} := 0$.
The symbol ``$\stackrel{D}{=}$'' 
will be used to denote equality in distribution
of two random variables, and ``$\stackrel{P}{\to}$''
will denote convergence in probability. Almost surely holding
statements will be abbreviated by ``a.s.'', and properties which
hold almost everywhere by ``a.e.''. The Lebesgue measure on
$\mathbb{R}$
will be denoted by $\lambda$. Throughout the paper, in order to
avoid 
trivialities, we will assume that $\xi$ and $\eta$ are different
from 
the zero process $t \mapsto 0$.

\section{Atoms of exponential \LL integrals} \label{s2}

Let $(\xi,\eta) = (\xi_t, \eta_t)_{t \geq 0}$ be a bivariate
L\'evy process.
Erickson and Maller~\cite{em1} characterised when the
exponential integral 
$I_t := \int_0^t e^{-\xi_{s-}} \, d\eta_s$, $t > 0$,
converges almost
surely to a finite random variable $I$ as $t\to\infty$. 
They showed that this happens if and only if 
\begin{equation} \label{eq-em1}
\lim_{t \to \infty} \xi_t = +\infty \; \; \mbox{a.s.}, \quad
\mbox{and}\quad
\int_{\mathbb{R} \setminus [-e,e]} 
\left( \frac{ \log |y| }{ A_\xi ( \log |y| ) } \right) \,
\Pi_\eta (dy) <
\infty. \end{equation} 
Here, the function $A_\xi$ is defined by 
$$A_\xi(y) := 1 + \int_1^y \Pi_\xi ((z,\infty)) \, dz, \quad y
\geq 1.$$ As
a byproduct of the proof, they obtained that $I_t$ converges
almost surely
to a finite random variable $I$ if and only if it converges in
distribution
to $I$, as $t\to \infty$. Observe that the convergence condition
\eqref{eq-em1} depends on the marginal distributions 
of $\xi$ and $\eta$ only, but
not on the
bivariate dependence structure of $\xi$ and $\eta$.

In this section we shall be interested in the question of
whether  
the limit random variable $I$ can have a  distribution with
atoms. A complete
characterisation of this will be given in Theorem~\ref{thm1}. A
similar
result for 
the characterisation of the existence of atoms for discrete time
perpetuities was obtained by Grincevi\v{c}ius~\cite{grin1},
Theorem~1. 
We will adapt his proof to show that $\int_0^\infty
e^{-\xi_{t-}} \,
d\eta_t$ has atoms if and only if it is  constant. This will be
a
consequence of the following lemma, which is formulated for
certain families
of 
random fixed point equations.

\begin{lemma} \label{lem-grin}
For every $t\geq 0$, let
$Q_t$, $M_t$ and
$\psi_t$ be random variables 
such that $M_t \not= 0$ a.s., and $\psi_t$ is independent
of $(Q_t, M_t)$.
Suppose $\psi$ is a random variable satisfying 
\[
\psi = Q_t + M_t \psi_t \quad {\rm for\ all }\;t\geq 0 ,
\]
and such that
\[
\psi  \=dr  \psi_t \quad {\rm for\ all}\; t\geq 0,
\]
and suppose further that
\[
Q_t  \topr  \psi \quad \mbox{as} \quad t \to \infty.
\]
Then $\psi$ has an atom if and only if it is 
a constant random variable.
\end{lemma}

\begin{proof}
We adapt  the proof of Theorem~1 of \cite{grin1}.
Suppose that $\psi$ has an atom at $a\in\mathbb{R}$, so that $$P
( \psi = a)
=: \beta > 0.$$ Then for all $\varepsilon \in (0,\beta)$ there exists some
$\delta >
0$ such that 
\begin{equation} \label{eq-1}
P (|\psi - a| < 2 \delta) < \beta  + \varepsilon. \end{equation}
Since $Q_t
\topr \psi$ as $t\to\infty$, there exists $t' = t'(\varepsilon)$
such that 
\begin{equation} \label{eq-2}
P(|\psi - Q_t| \geq \delta) = P( |M_t \psi_t| \geq \delta) <
\varepsilon\quad 
{\rm for\ all}\; t \geq t'.
\end{equation}
Then \eqref{eq-1} and \eqref{eq-2} imply that, for all $t \geq
t'$, \begin{equation} \label{eq-3} P( | Q_t - a| < \delta) \leq
P( |Q_t -
\psi| \geq \delta) + P(|\psi - a| 
< 2\delta) \leq \beta + 2 \varepsilon .
\end{equation}
Now observe that, for all $t\geq 0$,
\begin{eqnarray*}
\beta & = & P(\psi = a) \leq P( \psi=a, |\psi - Q_t| < \delta) 
+ P(|\psi - Q_t| \geq \delta) \\
& = & \int_{\mathbb{R}} P(Q_t + M_t s = a, |M_t s| < \delta) \,
dP(\psi_t
\leq s) + P(|\psi - Q_t| \geq \delta) \\ & = & \sum_{s\in D_t}
P(Q_t + M_t s
= a, |M_t s| < \delta) \, P(\psi_t = s) \;
+ P(|\psi - Q_t| \geq \delta).
\end{eqnarray*}
Here, the last equation follows from the fact that $P(Q_t + M_t
s =a)$ can
be positive for only a countable number of $s$, $s\in D_t$, say,
since the
number of atoms of any random variable is countable.

Since $\sum_{s\in D_t}P(\psi_t = s) \leq 1$ for all $s$, and
since $P(|\psi
- Q_t| \geq \delta) < \varepsilon$ for $t>t'$,  by
\eqref{eq-2}, it
follows that for such $t$ there is some $s_t \in \mathbb{R}$
such that
\begin{equation} \label{eq-4} \beta_t := P(Q_t + M_t s_t = a,
|M_t s_t| <
\delta) \ge \beta - \varepsilon. \end{equation} Observing that,
for all $t\geq 0$ $$\{ \psi = a\} \cup \{ Q_t + M_t s_t = a, |M_t s_t| < \delta
\} \subset
\{ |\psi - Q_t| \geq \delta \} \cup \{ |Q_t - a| < \delta \},$$
we obtain
for $t \geq t'$ that \begin{eqnarray*} \lefteqn{
P(|\psi - Q_t|
\geq \delta) + P(|Q_t - a| < \delta) }\\ & \geq & P(\psi = a) +
P(Q_t + M_t
s_t =a, |M_t s_t|  <\delta) - P(Q_t + M_t s_t = a, |M_t s_t|<
\delta, \psi =
a) \\ & = & \beta + \beta_t - \beta_t \, P(\psi_t = s_t).
\end{eqnarray*} We
used here that $P(M_t=0)=0$. {}From \eqref{eq-2} and
\eqref{eq-3} it now
follows that $$\beta_t \, P(\psi_t = s_t) \geq \beta + \beta_t -
\varepsilon
- 
(\beta + 2\varepsilon) = \beta_t - 3\varepsilon.$$
Using \eqref{eq-4} and the fact that $\psi \=dr \psi_t$,
we obtain $$P(\psi = s_t) = P(\psi_t = s_t) \geq 1 - \frac{3
\varepsilon}{\beta_t} > 
1 - \frac{3 \varepsilon}{\beta - \varepsilon}.$$
Letting $\varepsilon \to 0$ and observing that $P(\psi = a) >
0$, it follows
that $P(\psi = a) =1$. \end{proof}

As a consequence, we obtain:

\begin{theorem} \label{thm1}
Let $(\xi,\eta)$ be a bivariate L\'evy process such that $\xi_t$
converges
almost surely to $\infty$ as $t\to\infty$, 
and let $I_t := \int_0^t e^{-\xi_{s-}} \, d\eta_s.$
Denote the characteristic triplet of $(\xi,\eta)$  
by $({\gamma}, \Sigma, \Pi_{\xi, \eta})$, where ${\gamma}
= (\gamma_1, \gamma_2)$, and denote the upper diagonal element
of $\Sigma$
by $\sigma_\xi^2$. 
Then the following assertions are equivalent:
\begin{itemize}
\item[(i)] $I_t$ converges a.s. to a finite random variable $I$
as
$t\to\infty$, where $I$ has an atom. \item[(ii)] $I_t$ converges
a.s. to a
constant random variable as $t\to\infty$. 
\item[(iii)] $\exists$
$k\in\mathbb{R}\setminus \{ 0 \}$ such that $P\left(\int_0^t
e^{-\xi_{s-}} \, d\eta_s = k (1 -
e^{-\xi_t})\ {\rm for\ all}\ t > 0\right)=1$. 
\item[(iv)] $\exists$ $k\in
\mathbb{R}
\setminus \{ 0\}$ such that $e^{-\xi} = \mathcal{E} (-\eta/k)$,
i.e.
$e^{-\xi}$ is the stochastic 
exponential of $-\eta/k$.
\item[(v)] $\exists$ $k \in \mathbb{R} \setminus \{ 0\}$ such
that
$$\Sigma_{\xi, \eta} = \left( \begin{array}{cc} 1  & k \\ k &
k^{2} 
\end{array} \right) \sigma_\xi^2,$$
the L\'evy measure $\Pi_{\xi,\eta}$ of $(\xi,\eta)$ is
concentrated on $\{
(x, k (1-e^{-x}) ) : x \in \mathbb{R} \}$, \newline and 
\begin{equation} \label{claim}
\gamma_1 - k^{-1} \gamma_2 = \sigma_\xi^2 /2 + 
\int_{x^2 + k^2 (1-e^{-x})^2 \leq 1} (e^{-x} -1 + x) \,
\Pi_\xi(dx).
\end{equation} \end{itemize} \end{theorem}

\begin{proof}
To show the equivalence of (i) and (ii), suppose that 
$I$ exists a.s. as a finite random variable and define
\[
\psi := I = \int_0^\infty e^{-\xi_{s-}} \, d\eta_s, \quad
Q_t := I_t = \int_0^t e^{-\xi_{s-}} \, d\eta_s \quad \mbox{and}
\quad M_t :=
e^{-\xi_t}, \quad t \geq 0. \] Then \[ \psi \stackrel{D}{=}
\int_{t+}^\infty
e^{-(\xi_{s-} - \xi_t) } \, d (\eta_\cdot - \eta_t)_s =: \psi_t.
\]
So we have the setup of  Lemma~\ref{lem-grin}:
\begin{equation} \label{eq-MQ}
\psi = Q_t + M_t \psi_t,\ t\ge 0,
\end{equation}
 $Q_t$ converges in probability (in fact, a.s.) to $\psi$
as $t\to\infty$, and $\psi_t$ is independent of $(Q_t, M_t)$ for
all $t \geq
0$. We conclude from  Lemma~\ref{lem-grin} that $I = \psi$ is
finite a.s.
and has an atom if and only if it is constant, equivalently, if
(ii) holds.

Now suppose that (ii) holds and that the constant value of the
limit
variable is $k$. Then it follows from \eqref{eq-MQ} that, a.s.,
$$k=
\int_0^t
e^{-\xi_{s-}} \, d\eta_s + e^{-\xi_t} k, \ {\rm for\ each}\ t > 0,$$ hence 
\begin{equation} \label{eq-oscil}
\int_0^t e^{-\xi_{s-}} \, d\eta_s = k(1 - e^{-\xi_t}) \quad
{\rm for\ all}\; t > 0.
\end{equation} Observe that $k=0$ is impossible by uniqueness of
the 
solution to the
stochastic differential equation $d\int_0^t X_{s-} \, d\eta_s
=0$ (which
implies $e^{-\xi_s} = X_s = 0$, impossible). 
Since $Q_t$ and $e^{-\xi_t}$ are c\`adl\`ag functions,
\eqref{eq-oscil} holds on an event of probability 1.
This shows that (ii) implies (iii). 
The converse is clear, since $\lim_{t\to\infty}\xi_t =\infty$ 
a.s. by assumption.

Dividing \eqref{eq-oscil} by $-k$, we obtain
$e^{-\xi_t} = 1 + \int_0^t e^{-\xi_{s-}} d (-\eta_s/k)$, which
is just the
defining equation for $e^{-\xi} = \mathcal{E}(-\eta/k)$, see
Protter
\cite{protter}, p. 84, giving the equivalence of (iii) and (iv).

The equivalence of (iv) and (v) follows by straightforward
but messy calculations using the Dol\'eans-Dade formula and
the L\'evy-It\^o decomposition (for the calculation of
$\gamma$), and is
relegated to the  appendix. \end{proof}

\noindent{\bf Remarks.}\ (i)\
Under stronger assumptions,   Theorem \ref{thm1} may be
strengthened to conclude that $I$ has a density or is constant.
Suppose $(\xi,\eta)$ is a bivariate L\'evy process such that 
$\xi$ has no positive jumps and drifts to $\infty$, i.e.
$\lim_{t\to\infty}\xi_t=\infty$ a.s. Assume further that 
$\int_{\mathbb{R}\setminus [-e,e]}(\log |y|) \,\Pi_\eta(dy) <
\infty$.
  Then the condition
(\ref{eq-em1}) is fulfilled, and thus
$I:=\lim_{t\to\infty}\int_0^t e^{-\xi_{s-}}d\eta_s$ exists
and is finite
a.s. Applying the strong Markov property at the first passage
time
$T_x:=\inf\{t\geq0: \xi_t>x\} = \inf \{ t \geq 0: \xi_t = x \}$
(since $\xi$
has no positive
jumps)
yields the identity
$$I\,=\,\int_0^{T_x}e^{-\xi_{s-}}d\eta_s + e^{-x}I'$$
where $I'$ has the same distribution as $I$ and is independent of 
$\int_0^{T_x}e^{-\xi_{s-}}d\eta_s$. Thus $I$ is a
self-decomposable random
variable, 
and as a consequence its law is infinitely divisible
and unimodal and hence has a density, if it is not constant; see
Theorem
53.1, p.~404, in Sato \cite{sato}. Thus $I$ is continuous.  A
generalisation
of this result to the case of multivariate $\eta$ was recently
obtained
by 
Kondo et al.~\cite{KondoMaejimaSato}. 

(ii)\
As another important special case, suppose $\xi$ is a
Brownian motion
with a positive drift, and in addition that
$\int_{\mathbb{R}\setminus [-e,e]}(\log |y|) \,\Pi_\eta(dy) <
\infty$. Then
$I$ is finite a.s. {}From Condition (iii) of Theorem \ref{thm1}
we then see
that $\Delta \eta_t=0$, so the condition can hold only 
if $\eta_t$ is also a Brownian motion.
By Ito's lemma,  Condition (iii)  implies
$d\eta_t=k(d\xi_t-\sigma_\xi^2
dt/2)$, or, equivalently, $\eta_t=k(\xi_t-\sigma_\xi^2 t/2)$.
Similarly, if
$\eta$ is a Brownian motion, 
(iii) of Theorem \ref{thm1} can only hold 
if $\xi$ is a Brownian motion and the same relation is
satisfied. Thus we
can conclude that,  apart from this degenerate case, 
$\int_0^\infty e^{-B_s}d\eta_s$ and 
$\int_0^\infty e^{-\xi_s}dB_s$, when convergent a.s., have
continuous
distributions, for a Brownian motion $B_t$.

\section{Integrals with general $g$}\label{s3}

We now turn our  attention to the question of whether the
integral
$\int_0^\infty g(\xi_t)\, dt$ can have atoms, 
where $g$ is a more general deterministic function, and
$\xi = (\xi_t)_{t \geq 0}$ 
is a non-zero L\'evy process. 
To start with, we shall discuss some natural motivating
examples. Then we
shall present a few criteria that ensure the absence of atoms.
Finally, we
shall obtain by a different technique, which is a
variant of the stratification method,
a sufficient condition for the
absolute continuity of the integral.

\subsection{Some examples}
\begin{example} \label{ex-cp}
Let $(\xi_t)_{t \geq 0}$ be a compound Poisson process 
(with no drift) and 
$g: \mathbb{R} \to \mathbb{R}$ a deterministic function such
that
$g(0) \not=0$ and such that $\int_0^\infty g(\xi_t) \, dt$
is finite
almost surely. Then $\int_0^\infty g(\xi_t) \, dt$ has a
Lebesgue density.
\end{example}

\begin{proof}
Denote the time of the first jump of $\xi$  by $T_1$.
Recall that $\xi$ is always assumed nondegenerate, so $T_1$
is a nondegenerate exponential random variable. 
We can write $$\int_0^\infty g(\xi_t)
\, dt = g(0) T_1 + \int_0^\infty g(\xi_{T_1+t}) \, dt$$
(from which it is evident that the integral on the righthand
side converges
a.s.). Recall that the jump times in a compound Poisson process are
independent of the
jump sizes. By the strong Markov property of L\'evy processes
(see \cite{bertoin}, Prop. 6, p. 20), the process
$(\xi_{T_1+t})_{t \geq 0}$, and {\it a fortiori} the random variable  $\int_0^\infty g(\xi_{T_1+t}) \,
dt$, are independent of 
$T_1$. From
this follows the claim, since $g(0) T_1$ has a Lebesgue density
and hence
its sum with any independent random variable has also.
\end{proof}

The following example shows that this property does not carry
over to
compound Poisson processes with drift, at least not if the
support of $g$ is
compact.

\begin{example} \label{ex-comp-Poi-drift}
Let $\xi= (\xi_t)_{t \geq 0} = (at + Q_t)_{t \geq 0}$ be a
compound Poisson
process 
together with a deterministic drift $a \not= 0$, such that
$\lim_{t\to\infty}\xi_t={\rm sgn}(a)\infty$ a.s. Suppose that
$g$ is a
deterministic integrable Borel 
function with compact support. Then  $\int_0^\infty g(\xi_t)
\, dt$ is finite almost surely and its distribution has atoms.
\end{example}

\begin{proof}
Since $\xi$ drifts to $\pm \infty$ a.s., there is a random time
$\tau$ after
which $\xi_t \notin \supp  g$ for all $t$;  that is, if $\xi$
enters $\supp
g$ at all; if it doesn't, then $g(\xi_t)=0$ for all $t\ge 0$. In
either
case, 
$\int_\tau^\infty g(\xi_t) \, dt=0$,
and since $g$ is integrable and the number of jumps of $Q$ until
time $\tau$ is almost surely finite, it follows that 
$\int_0^\infty g(\xi_t) \, dt<\infty$ a.s.

Suppose now that $a>0$, so that  $\xi$ drifts to $+\infty$ a.s.,
and let
$r=\sup(\supp g)$. 
If $r\le 0$ there is a  positive probability that  $\xi$ does not
enter $\supp
g$, except, possibly, when $r=t=0$, and then  $g(\xi_t)=0$; in
either case,
$\int_0^\infty g(\xi_t) \, dt=0$ with positive probability,
giving an atom
at $0$. If $r>0$, let $T=r/a$. The event $A$ that the first jump
of $\xi$ 
occurs at or after  time $T$ has 
positive probability. On $A$, $\xi_t=at$ for all $0\le t\le T$.
Also, since
$\xi$ drifts to $+\infty$ a.s., on a subset of $A$ with positive
probability
$\xi$ does not re-enter $\supp g$ after time $T$. On this
subset, we  have
$\int_0^\infty g(\xi_t) \, dt = \int_0^{T} g(at) \, dt$, which
is constant.
Similarly if $a<0$. \end{proof}

Our third example relies on the following classical criterion
for the
continuity of infinitely divisible distributions (cf. Theorem
27.4, p.~175, in Sato
\cite{sato}), that we shall further use in the sequel.

\begin{lemma}\label{L1} Let $\mu$ be an infinitely divisible
distribution on
$\DoubleR$ with an infinite L\'evy measure, or with a 
non-zero Gaussian component. Then $\mu$ is
continuous.
\end{lemma}

If $\xi$ has infinite L\'evy measure, or no drift, Example
\ref{ex-comp-Poi-drift} may fail, as shown next:

\begin{example} \label{ex-s}
Suppose that $\xi$ is a subordinator with infinite L\'evy
measure, or is a
non-zero
subordinator with no drift. Then $\int_0^\infty 1_{[0,1]}
(\xi_t) \, dt$ is
finite a.s. and  has no atoms. \end{example}

\begin{proof}
Since $\xi_t$ drifts to $\infty$ a.s. it is clear that
$\int_0^\infty
1_{[0,1]}(\xi_t) \, dt$ is finite almost surely.
For $x > 0$ define
$$L_x := \inf \{ t > 0: \xi_t > x \}.$$
Then $\int_0^\infty 1_{[0,1]} (\xi_t) \, dt = L_1$, and for $a >
0$ we have
\begin{eqnarray*} \{ L_1 = a \} & = & \{ \inf \{ u: \xi_u > 1 \}
= a \} \\ &
= & \{ \xi_{a-\varepsilon } \leq 1 \; {\rm for\ all}\ \varepsilon > 0,
\quad \xi_{a
+ \varepsilon}
> 1 \; {\rm for\ all} \; \varepsilon > 0 \} \\
& \subseteq & \{ \xi_a = 1\} \cup \{ \Delta \xi_a> 0 \}.
\end{eqnarray*} A
\LL process is stochastically continuous so $P(\Delta
\xi_a>0)=0$. If  $\xi$
is a subordinator with infinite L\'evy measure, then $P(\xi_a =
1)=0$ by
Lemma \ref{L1}. Thus we get $ P(L_1=a)=0$. If  $\xi$ is a
subordinator with
no drift, then $\Delta \xi_{L_1}>0$ a.s. (\cite{bertoin}, p. 77)
(and this
includes the case of a compound Poisson), so again \[ P(L_1=a)=
P(L_1=a,\Delta \xi_{L_1}>0)\le P(\Delta \xi_a>0)=0. \]
\vskip-1cm
\end{proof}

\subsection{Some criteria for continuity} \label{s3-2}
We shall now present some fairly general criteria which ensure
the
continuity of the distribution of the integral
$\int_{0}^{\infty}g(\xi_{t})dt$ whenever the latter is finite
a.s.
and the L\'evy process $\xi$ is transient (see 
Bertoin~\cite{bertoin}, Section~I.4 or Sato~\cite{sato}, 
Section~35 for definitions and properties of transient and 
recurrent L\'evy processes).
\bigskip

\noindent{\bf Remarks.}\
(i)\ One might expect that the existence of $\int_0^\infty g(\xi_t)
\, dt$
already implies the transience of $\xi$. That this is not true in
general was shown by Erickson and Maller~\cite{em2}, Remark~(ii) after
Theorem~6. As a counterexample, we may take $\xi$ to be a compound
Poisson process with L\'evy measure $\Pi(dx)= \sqrt 2 \delta_{1} + \delta_{-\sqrt 2}$. Note that $\int x \Pi(dx)=0$, so
 $\xi$ is recurrent. Nonetheless $\xi$ never returns to $0$  after its first exit-time
and thus $0< \int_0^\infty {\bf 1}_{\{\xi_t =0\}}\, dt<\infty$ a.s.

(ii)\ Sufficient conditions
under which the existence of $\int_0^\infty g(\xi_t) \, dt$ implies the
transience of $\xi$ are mentioned in Remark (iii) after Theorem 6 of
\cite{em2}. One such sufficient condition is that there is some
non-empty open interval $J \subset \mathbb{R}$ such that 
$\inf \{ g(x) : x \in J\} > 0$.
\bigskip

We shall now turn to the question of atoms of $\int_0^\infty
g(\xi_t) \,
dt$.
For the next theorem, denote by $E^\circ$ the set of inner
points of 
a set $E$, by $\overline{E}$ its topological closure and by
$\partial E$ its
boundary.

\begin{theorem} \label{thm6}
Let $g:\mathbb{R} \to [0,\infty)$ be a deterministic Borel
function. Assume
that its support, $\supp g$, is compact, that $g > 0$ on $(\supp
g)^\circ$, and that $0 \in (\supp g)^\circ$. 
Write $\partial \supp g:=\supp g \setminus (\supp g)^\circ$
for the boundary of $\supp g$.
Let $\xi$ be a transient L\'evy process, and assume that
$I:=\int_0^\infty
g(\xi_t) \, dt$ is almost surely finite.
If either

\noindent {\rm (i)}  $\xi$
is of unbounded variation
and $\partial \supp g$ is finite,

\noindent or 

\noindent {\rm (ii)}  $\xi$ 
 is of bounded variation with 
zero drift and $\partial \supp g$ is at most countable,

\noindent then the distribution of $I$ has no atoms.
\end{theorem}

\begin{proof}
If $\xi$ is a compound  Poisson process without drift, the
result follows
from Example~\ref{ex-cp}, so we will assume that $\xi$ has
unbounded
variation, or is of bounded variation with zero drift such that
its L\'evy
measure is infinite, and that $g$ has the properties specified
in the
statement of the theorem.

Write 
$$I(x) := \int_0^x g(\xi_t) \, dt, \quad x \in (0,\infty]. $$
Then 
$x \mapsto I(x)$ is increasing and $I =
I(\infty)$ is finite a.s. by assumption, so $I(x)<\infty$ a.s.
for all $x\ge
0$. Plainly  $I(x)$ is a.s.
continuous at each
$x>0$.  Assume by way of contradiction that there is some $a
\geq 0$ such
that $P(I = a) > 0$, and proceed as follows.

Define 
$$T_s := \inf \{ u \geq 0: I(u) = s \},  \quad s \geq 0.$$
Since 
$\xi_t$ is adapted to the natural filtration
$\{{\cal F}_t\}_{t\ge 0}$ of $(\xi_t)_{t\geq 0}$,
so is $g(\xi_{\cdot})$ ($g$ is Borel), thus 
$\{T_s>u\} =\{ \int_0^u g(\xi_t)dt<s\} \in {\cal F}_u$,
because $I(\cdot)$ is adapted to $\{{\cal F}_t\}_{t\ge 0}$. Thus
$T_s$ is a
stopping time for each $s\ge 0$. Further, $T_s > 0$ for all $s >
0$. Since
$0 \in (\supp g)^\circ$, it is clear that $a\not= 0$. By
assumption, $\xi$ is
transient, so there is a finite random time $\sigma$ such that 
$\xi_t \not\in \supp g$ for all $t \geq \sigma$. Then
$I(\infty) = I(\sigma)$, and it follows that
$P\{T_{a}<\infty\} > 0$.

Define the stopping times $\tau_n := T_{a - 1/n} \wedge n$. Then
$(\tau_n)_{n \in \mathbb{N}}$ is strictly increasing to $T_{a}$,
showing
that $T_{a}$ is announceable; it follows that $t \mapsto \xi_t$
is continuous
at $t = T_a$ on $\{T_{a}<\infty\}$, see e.g. Bertoin
\cite{bertoin}, p. 21
or p. 39. 
Let $B=\{T_{a}<\infty,\  I(\infty) = a\}$. 
We restrict attention to $\omega\in
B$ from now on.
Since $T_a$ is the first time $I(\cdot)$ reaches $a$, for
every $\varepsilon > 0$ there must be
a subset $J_\varepsilon \subset (T_a - \varepsilon, T_a)$ of
positive
Lebesgue measure such that $g(\xi_t) > 0$ for all $t \in
J_\varepsilon$.
Thus $\xi_t\in \supp g$ for all $t\in J_\varepsilon$, and so
$\xi_{T_a}\in
\supp g$. Since we assume that $\partial \supp g := \supp g \setminus (\supp
g)^\circ$ is
countable, 
and that $\xi$ has infinite L\'evy measure or a non-zero
Gaussian component,
we have by Lemma \ref{L1} that $P(\xi_t \in \partial \supp g)=0$ for all $t>0$.
Consequently $$ E( \lambda \{ t \geq 0: \xi_t \in \partial \supp g \}) = 
\int_0^\infty P(\xi_t \in \partial \supp g)dt=0.
$$
It follows that there are times $t<T_a$ arbitrarily close to
$T_a$ with
$\xi_t$ in $(\supp g)^\circ$. By the continuity of $t \mapsto
\xi_t$ at $t =
T_{a}$, we then have $\xi_{T_{a}} \in \overline{ (\supp
g)^\circ}$ for
$\omega\in B'\subseteq B$, where $P(B')=P(B)>0$. Since $g > 0$
on $(\supp
g)^\circ$ it follows that 
$\xi_{T_{a}} \in \partial ((\supp g)^\circ)$
on the event $B'\subseteq\{ I(\infty) = a \} $;
for, if not, this would imply,
by an application of the Markov property, that
$I(t) > a$ for $t > T_{a}$, which is impossible. 

Now suppose (i), so that $\xi$ is of infinite variation. Then it follows from
Shtatland's
(1965) result (\cite{shtat}, see also Sato \cite{sato}, Thm
47.1, p. 351)
that 0 is regular for both $(-\infty,0)$ and $(0,\infty)$. Since
$\xi_{T_a}$ belongs to the finite set $\partial \supp g $,
there is an open interval $U \subset (\supp
g)^\circ$ which has $\xi_{T_a}$ either as left  or right end
point. In
either case, the regularity of 0 for $(0,\infty)$ and for
$(-\infty,0)$
implies that immediately after time $T_a$ there must be times
$t$ such that
$\xi_{t}$ is strictly less than $\xi_{T_a}$ and other times $t$
such that
$\xi_{t}$ is strictly greater than $\xi_{T_a}$. By the
continuity of $\xi$
at $T_a$, it follows that there must be times after $T_a$ such
that $\xi_t
\in U$. Consequently, there is some $\varepsilon =
\varepsilon(\omega)>0$
such that $\xi_{T_a
+ \varepsilon} \in (\supp g)^\circ$.
By the right-continuity of $\xi$ at $T_a +
\varepsilon$ it follows further 
that $I(\infty) > I(T_{a}) = a$ on  $B'$, where $P(B')>0$ and 
$B'\subseteq \{I(\infty) = a\}$, a contradiction.

Alternatively, suppose (ii), so that $\xi$ has finite variation and zero drift (and infinite
L\'evy measure).
Then it follows that $\xi$ almost surely does not hit single
points (by
Kesten's theorem \cite{kesten}; see \cite{bertoin}, p. 67). 
Thus, since $\partial ((\supp g)^\circ)\subseteq \supp g
\setminus (\supp
g)^\circ$ and the latter is at most countable, 
$\xi$ almost surely does not hit $\partial ((\supp g)^\circ)$. 
But on the set $B'$, where $P(B')>0$ and 
$B'\subseteq \{T_a<\infty,\ I(\infty) = a \}$,
we have $\xi_{T_{a}} \in \partial ((\supp g)^\circ)$,
contradicting $P(
I(\infty) = a) > 0$. \end{proof}

\noindent{\bf Remarks.}\ (i)\
The assumptions on the topological structure of 
$\{ x : g(x) > 0 \}$ in the previous theorem are easy to check.
That they cannot be completely relaxed can be seen from the
following example:
let
$g(x) = 1$ for all $x\in \mathbb{Q} \cap [-1,1]$ and $g(x) = 0$
otherwise,
then $\supp g = [-1,1]$, $(\supp g)^\circ = (-1,1)$, but $g>0$
on $(-1,1)$
does not hold. And in fact, it is easy to see that in that case
we have for
every L\'evy process of unbounded variation or infinite L\'evy
measure
that 
$$E \int_0^\infty g(\xi_t) \, dt = E \int_0^\infty 
1_{\mathbb{Q} \cap [-1,1]} (\xi_t) \, dt =
\int_0^\infty P( \xi_t \in \mathbb{Q} \cap [-1,1] ) \, dt = 0$$
by Lemma \ref{L1}, so that $\int_0^\infty g(\xi_t) \, dt = 0$
a.s.


(ii)\
Suppose $g$ is as in Theorem \ref{thm6}, 
and assume $\int_0^\infty g(x)dx<\infty$. Let
$\xi$ be a Brownian motion with non-zero drift. Then
$\int_0^\infty g(\xi_t)dt<\infty$ a.s. by Theorem 6 of
\cite{em2} and the
integral has a continuous distribution by Theorem \ref{thm6}.
\bigskip

Theorem \ref{thm6} allows a wide class of transient \LL
processes (we have
to exclude $\xi$ which are of bounded variation with nonzero
drift, by Ex.
\ref{ex-comp-Poi-drift}), 
but restricts us, essentially, to nonnegative $g$ which have
compact support. Another approach which 
combines excursion theory and Lemma \ref{L1} allows a much wider
class of
$g$ at the expense of placing restrictions on the local
behaviour of $\xi$.
Here is the first result in this vein. We refer  e.g. to
Chapters IV and V
in \cite{bertoin} for background on local time  and  excursion
theory for
L\'evy processes.

\begin{theorem}\label{T1} Let $g:\DoubleR\to[0,\infty)$ be a
measurable
function such that $g>0$ on some neighbourhood of $0$. Suppose
that $\xi$ is
a transient
L\'evy process such that $0$ is regular for itself, in the sense
that
$\inf\left\{t>0: \xi_t=0\right\}=0$ a.s., and that the integral
$I:=\int_{0}^{\infty}g(\xi_t)dt$ is finite a.s. Then the
distribution of $I$
has no atoms. \end{theorem}

\begin{proof} Thanks to Example 3.1, we may assume without
losing generality
that $\xi$ is not a compound Poisson. Then $0$ is an
instantaneous point, in
the sense that
$\inf\{t>0: \xi_t\neq0\}=0$ a.s. The assumption that
$\xi$ is transient implies that its last-passage time at $0$,
defined  by
$$\ell:=\sup\left\{t\geq0: \xi_t=0\right\},$$
is finite a.s.
Since the point $0$ is regular for itself, there exists a
continuous
nondecreasing local time process at level $0$ which we denote by
$L=(L_t,
t>0)$; we also introduce its right-continuous inverse
$$L^{-1}(t):=\inf\left\{s\geq0: L_s>t\right\}\,,\qquad t\geq0$$
with the
convention that $\inf\emptyset = \infty$. 
 The largest value of $L$, namely, $L_\infty$, is finite a.s.;
more
precisely, $L_{\infty}$ has an exponential distribution, and we
have
$L^{-1}(L_{\infty}-)=\ell$ and $L^{-1}(t)=\infty$ for every
$t\geq
L_{\infty}$ (\cite{bertoin}, Prop. 7 and Thm 8, pp. 113--115).
We denote the
set of discontinuity times of the inverse local time before
explosion by 
$${\cal D}:=\{t<L_{\infty}: L^{-1}(t-)<L^{-1}(t)\}$$ 
and then, following It\^o, we introduce for every $t\in{\cal D}$
the
excursion $\varepsilon(t)$ with finite lifetime $\zeta_t:=
L^{-1}(t)-L^{-1}(t-)$ by
$$\varepsilon_s(t):=\xi_{L^{-1}(t-)+s}\,,\qquad
0\leq s < \zeta_t\,.$$ It\^o's excursion theory shows that
conditionally on
$L_\infty$, the family of finite excursions $(\varepsilon(t),
t\in{\cal D})$
is distributed as the family of the atoms of a Poisson point
process with
intensity $L_{\infty}{\bf 1}_{\{\zeta<\infty\}}n$, where $n$
denotes the
It\^o measure of the excursions of the L\'evy process $\xi$
away from $0$,
and $\zeta$ the lifetime of a generic excursion (\cite{bertoin},
Thm 10, p.
118).

Since $\xi$ is not a compound Poisson process, the set of times
$t$ at which
$\xi_t=0$ has zero Lebesgue measure a.s., and we can express the
integral in
the form $I=A+B$ with \begin{equation}\label{eq??}
A:=\sum_{t\in{\cal
D}}\int_{L^{-1}(t-)}^{L^{-1}(t)}g(\xi_s)ds
=\sum_{t\in{\cal D}}\int_{0}^{\zeta_t}g(\varepsilon_s(t))ds\quad
\end{equation}
and
$$ B:=\int_{\ell}^{\infty}g(\xi_s)ds\,.$$
Excursion theory implies that $A$ and $B$ are
independent,
and hence we just need to check that $A$ has no atom.
Now, the conditional distribution of $A$ given $L_\infty$ is
infinitely
divisible, with L\'evy measure $\Lambda$ given by the image of
$L_{\infty}{\bf 1}_{\{\zeta<\infty\}}n$ under the map
$\varepsilon\rightarrow\int_{0}^{\zeta}g(\varepsilon_s)ds$.

The fact that $0$ is an instantaneous point implies that
the measure ${\bf 1}_{\{\zeta<\infty\}}n$ is infinite, and
further that the
excursions $\varepsilon(t)$ leave $0$ continuously for all
$t\in{\cal D}$
a.s. The assumption that $g>0$ on some neighbourhood of $0$ then
entails
that $\int_{0}^{\zeta_t}g(\varepsilon_s(t))ds>0$ for every
$t\in{\cal D}$.
Thus $\Lambda\{(0,\infty)\}=\infty$, and we conclude from Lemma
\ref{L1}
that the conditional distribution of $A$ given $L_\infty$ has no
atoms. 
It follows that $P(A=a) = E (P(A=a|L_\infty)) = 0$ for 
every $a > 0$, completing
the proof of our statement. 
\end{proof}

\noindent{\bf Remark.}\ 
See Bertoin \cite {bertoin}, Ch. V and Sato~\cite{sato},
Section~43, for discussions relevant to 
\LL processes for which $0$ is regular for itself.
\bigskip

An easy modification of the argument in Theorem \ref{T1}
yields the following
criterion in the special case when the L\'evy process has no
positive jumps.
This extends the result of Theorem \ref{thm6}  by allowing a
drift, as long
as there is no upward jump.

\begin{proposition}\label{P1}
Let $g:\DoubleR\to[0,\infty)$ be a measurable function
with $g>0$ on some neighbourhood of $0$. Suppose that
$\xi_t=at-\sigma_t$,
where $a>0$ and $\sigma$ is a subordinator with infinite L\'evy
measure
and no drift, and such that the
integral $I:=\int_{0}^{\infty}g(\xi_t)dt$ is finite a.s. 
Assume further that $a \not= E \sigma_1$, so that $\xi$ is
transient.
Then the distribution of $I$ has no atoms.
\end{proposition}

\noindent{\bf Remark.}\
We point out that in the case when $\xi$ is a L\'evy process
with no
positive jumps and infinite variation, then $0$ is 
regular for itself (\cite {bertoin}, Cor. 5, p. 192),
and thus Theorem \ref{T1} applies. Recall also Example 3.2 for
the case of
compound Poisson processes with drift.  Therefore our analysis
covers
entirely the situation when the L\'evy process has no positive
jumps and is
not the negative of a subordinator.

\begin{proof}  Introduce the supremum process
$\bar \xi_t:=\sup_{0\leq s \leq t}\xi_s$.
We shall use the fact that the reflected process $\bar \xi-\xi$
is Markovian
and that $\bar \xi$ can be viewed as its local time at $0$; see
Theorem
VII.1 in \cite{bertoin}, p. 189. The first-passage process
$T_x:=\inf\left\{t\geq0: \xi_t\geq x\right\}$ $(x\geq0)$
thus plays the role of the inverse local time.
It is well-known that $T_{\cdot}$ is a subordinator (killed at
some
independent exponential time when $\xi$ drifts to $-\infty$); 
more precisely, the hypothesis that $\xi_t=at-\sigma_t$ has
bounded
variation implies that the drift coefficient  of $T_{\cdot}$ is
$a^{-1}>0$.

Let us consider first the case when $\xi$ drifts to $\infty$, so
the
first-passage times $T_x$ are finite a.s. We write ${\cal D}$
for the set of
discontinuities of $T_{\cdot}$ and for every $x\in{\cal D}$, we
define the
excursion of the reflected L\'evy process away from 0 as
$$\varepsilon_s(x)\,=\,x-\xi_{T_{x-}+s}\,,\qquad 0\leq s <
\zeta_x:=T_x-T_{x-}\,.$$ 
According to excursion theory, the point measure
$$\sum_{x\in{\cal
D}}\delta_{(x,\varepsilon(x))}$$ is then a Poisson random
measure with
intensity $dx\otimes \bar n$, where $\bar n$ denotes the It\^o
measure of
the excursions of the reflected process $\bar \xi-\xi$ away from
$0$. Let
$b>0$ be such that $g>0$ on $[-b,b]$. We can express
$$\int_{0}^{\infty}g(\xi_s)ds \,=\,A+B+C$$ where
\begin{eqnarray*}
A\,&=&\,a^{-1}\int_{0}^{\infty}g(x)dx\,,\\
B\,&=&\,\sum_{x\in{\cal D},
x\leq b}\int_{T_{x-}}^{T_x} g(\xi_s)ds
\,=\,\sum_{x\in{\cal D},
x\leq b}\int_{0}^{\zeta_x} g(x-\varepsilon_s(x))ds\,,\\
C\,&=&\,\sum_{x\in{\cal D},
x>b}\int_{T_{x-}}^{T_x} g(\xi_s)ds
\,=\,\sum_{x\in{\cal D},
x> b}\int_{0}^{\zeta_x} g(x-\varepsilon_s(x))ds\,.
\end{eqnarray*}
The first term $A$ is deterministic, and $B$ and $C$ are
independent
infinitely divisible random variables (by the superposition
property of
Poisson random measures). More precisely, the  L\'evy measure of
$B$ is the
image of ${\bf 1}_{\{0\leq x \leq b\}}dx\otimes \bar n$ by the
map
$$(x,\varepsilon)\mapsto
\int_{0}^{\zeta}g(x-\varepsilon_s)ds\,.$$
Observe that the value of this map evaluated at any $x\in[0,b]$
and
excursion $\varepsilon$ is strictly positive (because excursions
return
continuously to $0$, as $\xi$ has no positive jumps). On the
other hand, the
assumption that the L\'evy measure of the subordinator
$\sigma_t=at-\xi_t$
is  infinite ensures that $0$ is an instantaneous point for the
reflected
process $\bar \xi-\xi$, and hence the It\^o measure $\bar n$ is
infinite. It
thus follows from Lemma \ref{L1} that the infinitely divisible
variable $B$
has no atom, which establishes our claim.

The argument in case $\xi$ drifts to
$-\infty$ is similar; the only difference is that the
excursion process is now stopped when an excursion with infinite
lifetime
arises. This occurs at time (in the local-time scale $\bar\xi$)
$\bar
\xi_{\infty}=\sup_{t\geq0}\xi_t$, where this variable has an
exponential
distribution. 
\end{proof}

\subsection{A criterion for absolute continuity} \label{s3-3}

Next we will investigate some different sufficient conditions, 
and some of them also ensure the existence of Lebesgue
densities. We will
work with more general integrals of the form $\int_0^\infty
g(\xi_t) \,
dY_t$ for a process $(Y_t)_{t \geq 0}$ of bounded variation, 
independent of the L\'evy process $\xi$.
The method will be a variant of the stratification method, 
by conditioning on almost every quantity apart from certain jump
times. Such an approach was also used by Nourdin and
Simon~\cite{ns}
for the study of absolute continuity of solutions to certain 
stochastic differential equations.

 We need the
following lemma, which concerns only deterministic functions.
Part (a) is just a rewriting of Theorem~4.2 in
Davydov et al.~\cite{dls}, and it is this part which will be
invoked when studying $\int_0^\infty g(\xi_t) \, dY_t$ for
$Y_t = t$.

\begin{lemma} \label{lem-an1}
Let $Y : [0,1] \to \mathbb{R}$ be a right-continuous
deterministic function of bounded variation. Let $f:[0,1] \to
\mathbb{R}$ be a deterministic Borel function such that 
\begin{equation} \label{eq-le2}
f \not= 0 \quad a.e.
\end{equation}
and such that the Lebesgue-Stieltjes integral $\int_0^1 f(t)
\, dY_t$ exists and is finite. Let 
$$H : (0,1] \to \mathbb{R}, \quad x \mapsto \int_{0+}^x f(t) \,
dY_t,$$ and
denote by $\mu := H(\lambda_{|{(0,1]}} )$ the image measure of
$\lambda$
under $H$. Then the following are sufficient conditions for 
(absolute) continuity of $\mu$:\\
$(a)$ Suppose the absolute continuous part of the measure
induced
by $Y$ on $[0,1]$ has a density which is different from zero
a.e.
Then $\mu$ is absolutely continuous.\\
$(b)$ Suppose that $Y$ is strictly increasing and that $f$ is in
almost every point $t \in [0,1]$ right- or left-continuous.
Then $\mu$ is continuous.
\end{lemma}

\begin{proof}
$(a)$ Denoting the density of the absolute continuous part
of $Y$ by $\phi$, it follows that $H$ is almost everywhere
differentiable with derivative $f \phi \not=0$ a.e., and the 
assertion follows from Theorem~4.2 in Davydov et al.~\cite{dls}.

(b) Suppose that $Y$ is strictly increasing and denote
$$K := \{ t \in (0,1) : f \, \mbox{is right- or left-continuous
in}\; t\}.$$
By assumption, $K$ has Lebesgue measure 1.
Using the right-/left-continuity, for every $t \in K$ such 
that $f(t) > 0$ there exists a unique maximal interval $J_+(t)
\subset (0,1)$
of positive length such that $t \in J_+(t)$ and  $f(y) > 0$ for all $y \in J_+(t)$.
By the axiom of choice
there exists a subfamily $K_+ \subset K$ such that
$(J_+(t): t \in K_+)$
are pairwise 
disjoint and their union covers $K \cap \{ t \in (0,1): f(t) > 0
\}$. Since
each of these intervals has positive length, there can only be 
countably many such intervals, so $K_+$ must be countable.

Similarly, we obtain a countable cover $(J_-(t): t \in K_-)$ of
$K \cap \{ t
\in (0,1) : f(t) < 0 \}$ with disjoint intervals. Now let 
$a \in \,\mbox{Range} (H)$. Then 
\begin{eqnarray*}
{H^{-1} (\{ a \})} &  \subset &
  \left(
\bigcup_{t \in K_+} (H^{-1}(\{ a \}) \cap J_+(t) )
 \right) 
\cup \left( \bigcup_{t \in K_-} (H^{-1}(\{ a \}) \cap J_-(t)
 ) \right)  \\
& & 
 \cup \left( [0,1] \setminus K \right) \cup \{ t \in [0,1] :
f(t)
 = 0 \} \cup \{ 0,1\}.
\end{eqnarray*}
Observing that 
$$\lambda \left( H^{-1} (\{ a \} ) \cap J_\pm(t) \right)
=  \lambda \left( (H_{| J_\pm(t)})^{-1} (\{ a \} )\right) = 0$$
since $H$ is strictly increasing (decreasing) on $J_+(t)$ ($J_-(t)$)
as a consequence
of $f > 0$ on $J_+(t)$ ($f < 0$ on $J_-(t)$) and strict increase of $Y$, it follows
that $\lambda( H^{-1} ( \{ a \} ) )= 0$, showing continuity of
$\mu$.
\end{proof}

We now come to the main result of this subsection. Note that
the case $Y_t=t$ falls under the case (i) considered in the
following
theorem, giving particularly simple conditions for absolute
continuity of $\int_0^\infty g(\xi_t) \, dt$. In particular,
part (b)
shows that if $\xi$ has infinite L\'evy measure and 
$g$ is strictly monotone on a neighbourhood of 0, 
then $\int_0^\infty g(\xi_t)\, dt$ is absolutely continuous.

\begin{theorem} \label{thm-3}
Let $\xi = (\xi_t)_{t \geq 0}$ be a transient
L\'evy process with non-zero L\'evy measure
$\Pi_\xi$. Let $Y = (Y_t)_{t \geq 0}$ be a stochastic process of
bounded variation on compacts which 
has c\`adl\`ag paths and which is
independent of $\xi$. Denote the density of the absolutely
continuous 
part of the measure induced by the paths $t \mapsto Y_t(\omega)$
by
$\phi_\omega$. Let $g : \mathbb{R} \to \mathbb{R}$ be a
deterministic
Borel function  
and suppose that the integral $$I := \int_{(0,\infty)} g(\xi_t)
\, dY_t
$$ exists almost surely and is finite.\\ 
$(a)$ {\em [general L\'evy
process]} 
Suppose that there are a compact interval 
$J \subset \mathbb{R} \setminus \{ 0 \}$
with $\Pi_\xi(J) > 0$ and some constant $t_0 > 0$ such that 
\begin{equation} \label{eq-T1}
\lambda ( \{  |t| \geq t_0: g(t) = g(t+z) \} ) = 0 \quad
{\rm for\ all}\; 
z \in J.
\end{equation}
{\rm Case (i):} If  $\lambda( \{ t \in [t_0,\infty): \phi(t) = 0
\} )= 0$
a.s.,
then $I$ is absolutely continuous.\\
{\rm Case (ii):} If $Y$ is strictly increasing on $[t_0,\infty)$
and $g$
has only countably many discontinuities, then $I$ does not have
atoms.\\
\\ $(b)$ {\em [infinite activity
L\'evy process]} Suppose the L\'evy measure $\Pi_\xi$ is
infinite. Suppose
further that there is $\varepsilon > 0$  such that
\begin{equation}
\label{eq-T2} \lambda ( \{ t \in (-\varepsilon,\varepsilon) :
g(t) = g(t+z)
\}
) = 0
\quad {\rm for\ all}\; z \in [-\varepsilon, \varepsilon].
\end{equation} 
{\rm Case (i):} 
If  $\lambda ( \{ t \in (0,\varepsilon) : \phi(t) = 0 \}) =0$
a.s.,
then $I$ is absolutely continuous.\\
{\rm Case (ii):} If $Y$ is strictly increasing on
$(0,\varepsilon)$ and $g$
has only countably many discontinuities, then $I$ does not have
atoms.
\end{theorem}

\begin{proof}
$(a)$ Let  $J$ be an interval such that \eqref{eq-T1} is
satisfied, and
define $$R_t := \sum_{0 < s \leq t, \Delta \xi_s \in J} \Delta
\xi_s, \quad
M_t := \xi_t - R_t, \quad t \geq 0.$$ Then $R = (R_t)_{t \geq
0}$ is a
compound Poisson process, independent of $M = (M_t)_{t \geq 0}$.
For $i\in\mathbb{N}$ 
denote by $T_i$ and $Z_i$ the
time and size of the $i^{\rm th}$ jump of $R$, respectively, 
and let $T_0 := 0$. Further, denote
\begin{eqnarray}
I_i & := & \int_{(T_{2i-2}, T_{2i}]} g(\xi_t) \, dY_t \nonumber
\\ & = &
\int_{(T_{2i-2}, T_{2i-1}]} \left( g\left(M_t + \sum_{j=1}^{2i-2} Z_j
\right)
- g\left( M_t + \sum_{j=1}^{2i-1} Z_j\right) \right) \, dY_t \nonumber \\ &
& + 
\int_{(T_{2i-2}, T_{2i}]} g\left(M_t + \sum_{j=1}^{2i-1} Z_j\right) \,
dY_t \nonumber
\\ & & + \left[ g( \xi_{T_{2i-1}} ) - g(\xi_{T_{2i-1}} -
Z_{2i-1} ) \right]
\Delta Y_{T_{2i-1}} \nonumber \\
& & + \left[ g( \xi_{T_{2i}} ) - g( \xi_{T_{2i}} - Z_{2i}
\right] \Delta Y_{T_{2i}}.\quad \label{eq-T2b}
\end{eqnarray}
We now condition on all random quantities present except the odd
numbered
$T_i$. Thus, 
for every Borel set $B \subset \mathbb{R}$, we write
$$P( I \in B) = E \, P \left( \sum_{i=1}^\infty I_i \in B \big|
Y, M, (T_{2j})_{j \in \mathbb{N}}, (Z_j)_{j \in \mathbb{N}}
\right).$$ To
show that $I$ has no atoms, it is hence sufficient to show that
\begin{equation} \label{eq-T3}  P \left( \sum_{i=1}^\infty I_i
\in B \big| 
Y, M, (T_{2j})_{j \in \mathbb{N}}, (Z_j)_{j \in \mathbb{N}}
\right) = 0 \quad {\rm a.s.}
\end{equation}
for every Borel set $B$ of the form $B = \{ a \}$ with $a \in
\mathbb{R}$.
Similarly, 
for showing that $I$ is absolutely continuous it is sufficient
to show 
that \eqref{eq-T3} holds for every Borel set $B$ of Lebesgue
measure $0$.
Observe that the $(I_i)_{i \in \mathbb{N}}$ are conditionally
independent
given $$V := (Y, M, (T_{2j})_{j \in \mathbb{N}}, (Z_j)_{j \in
\mathbb{N}}).$$ 
Thus the conditional probability that $I = \sum_{i=1}^\infty I_i\in B$
is
the convolution of the conditional probabilities that $I_i\in B$, $i
\in
\mathbb{N}$. Hence it suffices to 
show that there is some random integer $i_0 \in \mathbb{N}$ such
that almost
surely, the conditional distribution of $I_{i_0}$ given $V$ is
absolutely continuous (case (i)) or has no atoms (case (ii)),
respectively. 
 
Define the integer $i_0$ as the first index $i$ such that
\begin{equation}
\label{eq-T4} \min \left\{  \inf_{ t \in (T_{2 i - 2}, T_{2i} ]}
\{ |M_t + \sum_{j=1}^{2 i -2} Z_j | \} \, , T_{2i-2} \right\}
\geq t_0,
\end{equation}
with $t_0$ as in \eqref{eq-T1}. Since $\xi$ is transient
$i_0$ is almost surely finite.
As a function of $V$,  $i_0$ is constant under the conditioning
by $V$. The
right hand side of  \eqref{eq-T2b} is comprised of  four
summands. The
second and fourth summands are constant given $V$. The third
summand is
still random, after conditioning, since $T_{2i-1}$ enters in
$\Delta Y$; but
here $R$ and $Y$ 
are independent, so that the third summand equals 0 a.s.
Thus it is sufficient to show that, given $V$, the first summand,
evaluated at $i_0$, namely $$\widetilde{I}_{i_0} := \int_{
(T_{2i_0 - 2},
T_{2 i_0 - 1}] } \left( g\left(M_t + \sum_{j=1}^{2 i_0 -2} Z_j\right) - 
g\left(M_t + \sum_{j=1}^{2 i_0 - 1} Z_j\right) \right)dY_t,$$
is almost surely absolutely continuous (case (i))
or has no atoms (case (ii)). Define
the functions 
$f = f_V: [T_{2i_0 - 2}, T_{2i_0}] \to \mathbb{R}$ and
$H = H_V : (T_{2_{i_0} - 2}, T_{2i_0}] \to \mathbb{R}$ by
\begin{eqnarray*}
f(t) & = & g\left(M_t + \sum_{j=1}^{2i_0 - 2} Z_j \right) - \
g\left(M_t + \sum_{j=1}^{2 i_0 - 1} Z_j \right),\\
H(x)  &:= & \int_{(T_{2i_0 - 2}, x]} f(t) \, dY_t.
\end{eqnarray*} Observing
that $T_{2i_0 - 1}$ is uniformly distributed on 
$(T_{2i_0 - 2}, T_{2i_0})$ given $V$, it follows from
Fubini's theorem that for any Borel set $B \subset \mathbb{R}$ 
\begin{eqnarray*}
P (\widetilde{I}_{i_0} \in B | V) & = & 
E ( \mathbf{1}_{\{ H(T_{2{i_0} -1} ) \in B \} } | V) =
\int_{(T_{2 i_0 - 2}, T_{2 i_0} )} \mathbf{1}_{
\{ H(x) \in B \} } \, \frac{ dx}{T_{2 i_0}
- T_{2 i_0 - 2} } \\
& = & \frac{ \lambda (H^{-1} (B))}{ T_{2 i_0} - T_{2 i_0 - 2} }.
\end{eqnarray*} We shall apply Lemma~\ref{lem-an1} to show that
$\widetilde{I}_{i_0}$ given $V$ is absolutely continuous or has
no atoms,
respectively. For this, observe that \eqref{eq-le2} is satisfied
because of
\eqref{eq-T1} and \eqref{eq-T4}, and note that $z:=Z_{2i_0-1}
\in J$, since
all the jumps of $R$ are in the interval $J$.
In case (i) this then
gives absolute continuity of $\widetilde{I}_{i_0}$ conditional on $V$ by
Lemma~\ref{lem-an1}~(a) and hence of the distribution of $I$.
Now concentrate
on case (ii), when $Y$ is strictly increasing on $[t_0,\infty)$
and 
$g$ has only countably many discontinuities. Denote
this set of 
discontinuities of $g$ by $F$. By assumption, $F$ is countable.
This then implies that almost surely, the function $f$ is almost
everywhere 
right-continuous.  
For by the a.s. right-continuity of the paths of 
L\'evy processes, $f$ can happen to be non-right-continuous at a
point 
$t$ only if $\xi_t^{(1)} := M_t 
+ \sum_{j=1}^{2 i_0 - 1} Z_j \in F$ or $\xi_t^{(2)} := M_t +
\sum_{j=1}^{2 i_0 - 2}
Z_j \in F$. But 
$$E ( \lambda \{ t \geq 0 : \xi_t^{(1)} \in F \; \mbox{or}\; 
\xi_t^{(2)} \in F \} ) =
\int_0^\infty P (\xi_t^{(1)} \in F \; \mbox{or} \; \xi_t^{(2)}
\in F
) \, dt,$$
and by Lemma \ref{L1}
the last integral is zero if $\xi$ has infinite L\'evy measure,
so that almost surely, $f$
is almost everywhere right-continuous if $\Pi_\xi$ is infinite.
If $\xi$ has finite L\'evy measure, then $f$ is trivially
almost everywhere right-continuous. So we see that in case (ii)
our assumptions imply the
conditions 
of Lemma~\ref{lem-an1}~(b), which then gives the claim.

$(b)$ The proof is similar to the proof of $(a)$: 
for $0 < \delta < \varepsilon / 2$, let
$$R_t^{(\delta)} := \sum_{ |\Delta \xi_s| \in [\delta,
\varepsilon/2] } 
\Delta \xi_s, \quad M_t^{(\delta)} := \xi_t -
R_t^{(\delta)},\quad t \geq
0,$$ and denote the time and size of the $i^{\rm th}$ jump of
$R^{(\delta)}
= (R_t^{(\delta)})_{t \geq 0}$ by 
$T_i^{(\delta)}$ and $Z_i^{(\delta)}$, respectively.
Further, define the set $\Omega_\delta$ by
$$\Omega_\delta := \{ T_2^{(\delta)} \leq \varepsilon, \sup_{0
\leq t <
T_2^{(\delta)} }
| M_t^{(\delta)} | \leq \varepsilon/2 \}.$$
Let $P_\delta (\cdot ) := P( \cdot | \Omega_\delta)$, and denote
expectation
with respect to $P_\delta$ by $E_\delta$. 
Since $P(\Omega_\delta) \to 1$ as $\delta \downarrow 0$ because
the L\'evy measure of $\xi$ is infinite,      
it is sufficient to show that, given $\delta > 0$, we have
$P_\delta (B) =
0$ for  all Borel sets $B$ such that $\lambda(B) = 0$
(case (i)), or such that 
$B= \{a\}$, 
$a \in \mathbb{R}$ (case (ii)),  respectively. Let 
$$V_\delta := (Y, M^{(\delta)}, (T_j^{(\delta)})_{j \geq 2},
(Z_j^{(\delta)})_{j \in \mathbb{N}}.$$ Then we can write
$$P_\delta (I \in
B) = E_\delta P_\delta (I \in B | V_\delta),$$ and it suffices
to show that
$P_\delta (I \in B| V_\delta) = 0$ a.s. for the sets $B$ under
consideration. But, conditional on $V_\delta$, $I$ almost surely
differs
from 
$$\widetilde{I}_1 := \int_{(0, T_2^{(\delta)}]} \left( 
g\left(M_t^{(\delta)} \right) - g\left( M_t^{(\delta)} + Z_1^{(\delta)} \right)
\right) \, dY_t$$
only by a constant. It then follows in complete analogy to the
proof 
of $(a)$ that under $P_\delta$, $\widetilde{I}_1$ given
$V_\delta$ has no
atoms or is absolutely continuous, respectively, which then
transfers to $I$
under $P_\delta$ and hence to $I$ under $P$. \end{proof}

\noindent {\bf Remarks.}\ (i)\
The preceding proof has shown that the
independence assumption on $\xi$ and $Y$ can be weakened.
Indeed,  we need only assume 
that the processes $(R_t)_{t \geq 0}$ and $Y$ are independent.\\
(ii)\ 
In addition to the assumptions of Theorem \ref{thm-3},
assume that $g$ is continuous. Then almost surely,
$g(\xi_{t-}) = g(\xi_{t})_-$ exist for all $t > 0$, and
the assertions of  Theorem~\ref{thm-3} 
remain true for integrals of the form
$$\int_{(0, \infty)} g(\xi_t)_- \, dY_t.$$
This follows in complete analogy to the proof of Theorem
\ref{thm-3}.\\
(iii)\ Similar statements as in Theorem~\ref{thm-3} can be made
for integrals of the form
$\int_0^\infty (g(\xi_t + \psi(t)) \, dt$,
where $\psi$ is some deterministic function behaving nicely.
We omit the details.

\section*{Appendix}

\noindent
{\bf Proof of the equivalence of (iv) and (v) in
Theorem~\ref{thm1}.}
Assume (iv), and observe that by the 
Dol\'eans-Dade formula (e.g.~\cite{protter}, p.~84), $e^{-\xi} =
\mathcal{E} 
(-\eta/k)$, where $k \not= 0$, if and only if $\Pi_\eta ( \{ y \in \mathbb{R} :
k^{-1} y \geq 1
\} ) = 0$ and  $\xi_t = X_t$, where \begin{equation}
\label{dade} X_t :=
k^{-1} \eta_t + k^{-2} \sigma_\eta^2 t/2 - \sum_{0 \leq s \leq
t} 
\left( \log (1 - k^{-1} \Delta \eta_s ) + k^{-1} \Delta \eta_s
\right),
\quad t \geq 0. \end{equation} Now $(X, \eta)$ is a bivariate
L\'evy
process, whose 
Gaussian covariance matrix is given by
$\Sigma_{X, \eta} = \left(
\begin{array}{cc}
1  & k \\
k & k^{2} 
\end{array} \right) \sigma_X^2$.
Further, \eqref{dade} implies $\Delta X_t = -\log (1 - k^{-1}
\Delta
\eta_t)$, showing that the L\'evy measure $\Pi_{X,\eta}$ of 
$(X, \eta)$ is concentrated on $\{ (x, k (1-e^{-x}) ): x\in
\mathbb{R} \}$.

Conversely, if $(Y,\eta)$ is a bivariate L\'evy  process with
Gaussian covariance matrix  given by $\Sigma_{Y,\eta}=\Sigma_{X, \eta}$, whose
L\'evy measure is concentrated on $\{ (x, k (1-e^{-x}) ): x\in
\mathbb{R} \}$, then  $\Delta Y_t = -\log (1 - k^{-1}\Delta
\eta_t)$,  and it follows that there is some $c \in \mathbb{R}$
such that $Y_t = X_t + ct$, so that $e^{-Y_t + ct}=
(\mathcal{E}(-\eta/k))_t$.
Hence we have established the
equivalence of (iv) and (v)
in Theorem \ref{thm1}, subject to relating  $\gamma_1$ and
$\gamma_2$ as in
\eqref{claim}.

To do this, let $X_t$ as in \eqref{dade}
and  use the L\'evy--It\^{o} decomposition. Define
$$
\left( { X_t^{(1)} \atop \eta_t^{(1)} } \right) :=
\lim_{\varepsilon
\downarrow 0} \left( 
\sum_{0 < s \leq t \atop (\Delta X_s)^2 + (\Delta \eta_s)^2 >
\varepsilon^2}
\left( { \Delta X_s \atop \Delta \eta_s } \right) - t \int
\int_{x_1^2 +
x_2^2 \in (\varepsilon^2, 1]} 
\left( { x_1 \atop x_2 } \right) \Pi_{X, \eta} (d (x_1, x_2))
\right)$$ and
$(X_t^{(2)}, \eta_t^{(2)})' := (X_t, \eta_t)' - (X_t^{(1)},
\eta_t^{(1)})'$ where the limit is a.s. as $\varepsilon \downarrow 0$.
(Note that the expression in big brackets on the right is not
precisely the
compensated sum of jumps.) Then $(X_t^{(2)},
\eta_t^{(2)})_{t\geq 0}'$ is a
L\'evy process with characteristic triplet $(\gamma, \Sigma,
0)$, so has the
form $(X_t^{(2)}, \eta_t^{(2)})' = (\gamma_1 t, \gamma_2 t)' +
\vec{B}_t$,
$t \geq 0$, where $(\vec{B}_t)_{t \geq 0}$ is a Brownian motion
in
$\mathbb{R}^2$. {}From this follows that 
\begin{equation} \label{eq-compare}
X_t^{(2)} - k^{-1} \eta_t^{(2)} = (\gamma_1 - k^{-1} \gamma_2) t
+ \widetilde{B}_t, \quad
t \geq 0,
\end{equation}
for some Brownian motion $(\widetilde{B}_t)_{t \geq 0}$ in
$\mathbb{R}^1$.
We wish to determine $\gamma_1 - k^{-1} \gamma_2$. To do this,
observe that
from \eqref{dade} and $\sigma_X^2 = k^{-2} \sigma_\eta^2$, we
have
\begin{eqnarray*}&&(X_t - X_t^{(1)}) - k^{-1} (\eta_t -
\eta_t^{(1)})\\ &&=\sigma_X^2 t/2 + \sum_{0 \leq s \leq t}
(\Delta X_s
- k^{-1} \Delta \eta_s) \\ & & - \lim_{\varepsilon \downarrow 0}
\left( 
\sum_{0 < s \leq t \atop (\Delta X_s)^2 + (\Delta \eta_s)^2 >
\varepsilon^2}
(\Delta X_s - k^{-1} \Delta \eta_s) - t
\int \int_{x_1^2 + x_2^2 \in (\varepsilon^2, 1]} (x_1 - k^{-1}
x_2) 
 \, \Pi_{X, \eta} (d (x_1, x_2)) \right).
 \end{eqnarray*}
 Noting that $k^{-1} \Delta \eta_s = 1- e^{-\Delta X_s}$ and
that 
 $\sum_{0 < s \leq t} (\Delta X_s - 1 + e^{-\Delta X_s})$
converges
absolutely, we obtain, letting $\varepsilon \downarrow 0$, that
\begin{eqnarray*}  X_t^{(2)} - k^{-1} \eta_t^{(2)} & = & 
 \sigma_X^2 t/2 + t \int\int_{x_1^2 + x_2^2 \leq 1} (x_1 -
k^{-1} x_2) 
 \Pi_{X, \eta} (d (x_1, x_2)) \\
 & = & \sigma_X^2 t/2 + t \int_{x^2 + k^2 (1-e^{-x})^2 \leq 1}
 (x-1+e^{-x}) \, \Pi_X (dx).
 \end{eqnarray*}
 Comparing this with \eqref{eq-compare}  gives \eqref{claim}.
$\halmos$

\subsection*{Acknowledgements}
This research was carried out while JB and AL were visiting
the Centre for Mathematical Analysis and the School of Finance
\& Applied Statistics at ANU in Canberra. They take pleasure
in thanking both for their hospitality. AL gratefully
acknowledges
financial support by the German Science Foundation (Deutsche 
Forschungsgemeinschaft),  research grant number Li 1026/2-1,
and RM's research was partially supported by ARC grant DP00664603.

\end{document}